\documentclass[12pt]{amsart}

\def \Z {{\mathbb Z}}
\def \N {{\mathbb N}}
\usepackage{amsmath,amssymb}
\newtheorem{theorem}{Theorem}[section]

\newtheorem{proposition}[theorem]{Proposition}
\newtheorem{lemma}[theorem]{Lemma}

\newtheorem{definition}[theorem]{Definition}

\begin{document}
\title{Star-free geodesic languages for groups}

\author[S.~Hermiller]{Susan Hermiller$\!\,^1$}
\address{Dept. of Mathematics\\
        University of Nebraska\\
         Lincoln NE 68588-0130  USA}
\email{smh@math.unl.edu}

\author[Derek F.~Holt]{Derek F.~Holt}
\address{Mathematics Institute\\
      University of Warwick  \\
       Coventry CV4 7AL, UK  }
\email{dfh@maths.warwick.ac.uk}

\author[Sarah Rees]{Sarah Rees}
\address{School of Mathematics and Statistics\\
       University of Newcastle \\
       Newcastle NE1 7RU, UK  }
\email{Sarah.Rees@ncl.ac.uk}

\date{\today}

\maketitle

\footnotetext[1]{Supported under NSF grant no.\ DMS-0071037}

\begin{abstract} 
In this article we show
that every group with a finite presentation
satisfying one or both of the small cancellation conditions
$C'(1/6)$ and $C'(1/4)-T(4)$ has the property that the
set of all geodesics (over the same generating set)
is a star-free regular language.  Star-free regularity
of the geodesic set is shown to be dependent on
the generating set chosen, even for free groups.  
We also show that the class of groups 
whose geodesic sets are star-free with respect
to some generating set is closed under
taking graph (and hence free and direct) products,
and includes all virtually abelian groups. 
\end{abstract}

\section{Introduction}\label{intro}

There are many classes of finitely presented groups that have been
studied via sets of geodesics that are regular languages 
(that is, sets defined
by finite state automata). 
Various examples are known of groups for which the
set of all geodesics is a regular set. 
Word hyperbolic groups with any
finite generating set are very natural
examples; for these the sets of geodesics
satisfy additional (``fellow traveller'') properties which make the
groups automatic~\cite[Theorem 3.4.5]{ECHLPT}.
Other examples include finitely generated abelian groups 
with any finite generating
set~\cite[Propositions 4.1 and 4.4]{NeumannShapiro}; and
with appropriate generating sets, virtually abelian
groups, geometrically finite hyperbolic groups
~\cite[Theorem 4.3]{NeumannShapiro},
Coxeter groups (using the  standard generators)~\cite{Howlett},
Artin groups of finite type and, more generally, Garside groups
~\cite{CharneyMeier} 
(and hence torus knot groups).

To date very little connection has been made between the properties of the
regular language that can be associated with a group in this way and the
geometric or algebraic properties of the group itself.
In this paper we consider groups whose sets of geodesics satisfy 
the more restrictive language theoretic property of star-free
regularity.

The set of regular languages over a finite alphabet $A$ 
is by definition the closure
of the set of finite subsets of $A$ under the operations of union,
concatenation, and Kleene closure~\cite[Section 3.1]{HMU}.
(The Kleene closure $X^*$ of a set $X$ is defined to be the
union $\cup_{i=0}^\infty X^i$ of concatenations of copies of $X$ with itself.)
By using the fact that the regular languages are precisely those that
can be accepted by a finite state automaton, it can be proved that
this set is closed under many more operations including
complementation and intersection~\cite[Section 4.2]{HMU}.
The set of {\it star-free languages} over $A$ is defined to be
the closure of
the finite subsets of $A$ under concatenation and the three Boolean
operations of union,
intersection, and complementation, but not under Kleene 
closure~\cite[Chapter 4 Definition 2.1]{Pin}.

The star-free languages form a natural low complexity subset of the
regular languages.  An indicator of the fundamental role that they play
in formal language theory is the surprisingly large variety of
conditions on a regular language that turn out to be equivalent to
that language being star-free.
The book by McNaughton and Papert~\cite{MP} is devoted entirely to
this topic. One such condition that we shall make use of in this paper
is the result of Sch\"utzenberger that a regular language is star-free if and
only if its syntactic monoid is aperiodic.
Other examples studied in~\cite{MP} are the class $LF$
of languages represented by nerve nets that are buzzer-free and almost
loop-free, and the class $FOL$ of languages defined by a sentence in first
order logic.  A regular language not containing the empty string is star-free
if and only if it lies in $LF$ or, equivalently, in $FOL$.
There are also relationships between star-free languages
and various types of Boolean circuits; these are discussed
in detail in the book by Straubing~\cite{Straubing}.

Note that $A^*=\emptyset^c$ is star-free.  
The set of star-free languages
properly contains the set of all {\em locally testable} languages.
A subset of $A^+ :=A^*\setminus \{\epsilon\}$ (where $\epsilon$ is the
empty word) is locally testable if it is
defined by a regular expression which combines terms of the form
$A^*u$, $vA^*$ and $A^*wA^*$, for non-empty finite strings $u,v,w$,
using the three Boolean operations~\cite[Chapter 5 Theorem 2.1]{Pin}. 
A natural example of a locally testable language is provided by the set of all
non-trivial geodesics in a free group in its natural presentation.
Further examples of star-free languages are provided by the {\em piecewise
testable} languages; a subset of $A^*$ is piecewise testable if
it is defined by a regular expression combining terms
of the form $A^*a_1A^*a_2 \cdots A^*a_kA^*$, where $k \ge 0$
and each $a_i \in A$, using the three Boolean 
operations~\cite[Chapter 4 Proposition 1.1]{Pin}.
A natural example of a piecewise testable
language is provided by the set of geodesics in a free abelian group
over a natural generating set.

Margolis and Rhodes~\cite{MargolisRhodes} conjectured that 
the set of geodesics in any word hyperbolic group
(with respect to any generating set) is star-free.
This conjecture was motivated by an interpretation of
van Kampen diagrams in terms of Boolean circuits, and
utilization of relationships between properties of
circuits and star-free languages mentioned above.

We show that the conjecture as stated is false in Section~\ref{dependence},
where we describe a 6-generator presentation of the free group of
rank 4 for which the set of geodesics is not star-free.
This demonstrates that the property of having star-free geodesics
must be dependent on the choice of generating set.

By contrast, our main theorem,
in Section~\ref{sc}, states that groups defined by a presentation satisfying
either one of the small cancellation conditions $C'(1/6)$ or $C'(1/4)-T(4)$
(which imply but are not a consequence of hyperbolicity)
have star-free sets of geodesics with respect to the generating set of
that presentation.  
Our proof relies on the thinness of van Kampen
diagrams for these presentations, which was used
in ~\cite{Strebel} to show that these groups are word hyperbolic.
The question of whether or not any word hyperbolic group must have some
generating set with respect to which the geodesics are star-free
remains open.

In Section~\ref{grprod} we consider closure properties of
the set of groups that have star-free sets of geodesics for some
generating set.
We show that
this set is closed under taking direct, free, and graph products.
Our proof follows the strategy of the proof in \cite{LMW}
showing that the set of groups with regular languages of
geodesics (for some generating set) is closed under graph products.

It is proved in~\cite[Propositions 4.1 and 4.4]{NeumannShapiro} that
every finitely generated abelian group has a
regular set of geodesics with respect to any finite generating set, and that
every finitely generated
virtually abelian group $G$ has a regular set of geodesics with respect to
some finite generating set.
(An example of J.W. Cannon,
described after Theorem
4.3 of~\cite{NeumannShapiro}, shows that 
the set of geodesics need not be regular
for every finite generating set of 
a virtually abelian group.)
In Section~\ref{va} we strengthen these results to show
that all finitely generated abelian groups have piecewise testable
(and hence star-free) sets of geodesics with respect to 
any finite generating set,
and that all finitely generated
virtually abelian groups also have piecewise testable sets of geodesics with
respect to certain finite generating sets.

\section{Some technicalities and basic results}

Let $A$ be a finite alphabet and
let $M$ be a deterministic finite state
automaton over $A$.  We assume that all such automata in
this paper are complete; that is, for every state and letter in
$A$, there is a corresponding transition.
For  a state $\sigma$ of $M$ and word $u \in A^*$, 
define $\sigma^u$ to be the
state of $M$ reached by reading $u$ from state $\sigma$.
The {\em transition monoid} associated with $M$ is defined
to be the monoid of functions between the states of $M$ induced by
the transitions of $M$.  Since the automaton is finite, this
transition monoid is a finite monoid.

For a regular language $L$ of $A^*$, there is a
minimal finite state automaton $M_L$ accepting the
language $L$, which is unique up to the naming of the states~\cite{HMU}.
In this minimal automaton, no two states have the same ``future'';
that is, the two sets of words labelling transitions from two distinct states
to accept states must be distinct.
The {\em syntactic congruence} of $L$ is the congruence
$\sim_L$, which relates two words $u$ and $v$ provided, 
for all words $x$ and $y$,
$xuy \in L$ if and only if $xvy \in L$.  The {\em syntactic monoid}
associated with $L$ is the quotient monoid $A^*/\sim_L$.
The images $[u]$ and $[v]$ of $u$ and $v$ in the
syntactic monoid $A^*/\sim_L$ are equal if and only if
$\sigma^u=\sigma^v$ for all states $\sigma$ of $M_L$.
Thus the syntactic monoid of $L$ is also the transition
monoid of $M_L$.

A monoid is said to be {\em aperiodic} if it satisfies a rule $x^N=x^{N+1}$
for some $N \in \N$.
By a result of Sch\"utzenberger \cite[Chapter 4 Theorem 2.1]{Pin},
a regular language $L$ over $A^*$ is star-free
if and only if its syntactic monoid is aperiodic.
The following proposition yields an alternative interpretation of
star-free, which allows an easy algorithmic check (which we used
to check examples).

\begin{proposition}\label{star-free1}  Let $L$ be a regular language
over $A$.
The language $L$ is not star-free if and only if, for any
$N \in \N$, there exist $n > N$ and
words $u,v,w \in A^*$ with one of the words
$uv^nw$ and $uv^{n+1}w$  in $L$ and
the other not in $L$.  Moreover, if $L$ is not star-free,
then there exist
fixed words $u,v,w \in A^*$ such that, for each $N \in \N$, there exists
$n > N$ with $uv^nw \in L$ but $uv^{n+1}w \not\in L$. 
\end{proposition}

\begin{proof}
If $L$ is star-free, Sch\"utzenberger's Theorem says that the
syntactic monoid is aperiodic, so there is an $N \in \N$ with
$x^N=x^{N+1}$ for all $x$ in the monoid $A^*/\sim_L$.  
For any $u,v,w \in A^*$,
the image $[v]$ of $v$ in $A^*/\sim_L$
satisfies $[v]^N=[v]^{N+1}$, and so
$\sigma^{v^N}=\sigma^{v^{N+1}}$ for all states $\sigma$
of the minimal automaton $M_L$ accepting $L$, including
the state $\sigma=\sigma_0^u$, where $\sigma_0$ is the start state.
Thus for all $n>N$, the word $uv^nw$ will be accepted by $M_L$ if and only if
$uv^{n+1}w$ is accepted.  The proof of the converse is similar.

Now suppose that 
$L$ is not star-free, and let $m$ be the number of
states of the automaton $M_L$.
The first part of the proposition shows that
there are words $u,v,w \in A^*$ and an integer $i>m$
with
one of the words $uv^{i}w$ or $uv^{i+1}w$ in $L$
and the other in $A^* \setminus L$. 
Since $i>m$ there must be two integers $j_1,j_2$ with $0 \leq j_1<j_2<i$
such that 
$\sigma_0^{uv^{j_1}}=\sigma_0^{uv^{j_2}}$. Let $k=j_2-j_1$.
Then for all natural numbers $l$, 
$\sigma_0^{uv^{i}w}=\sigma_0^{uv^{i+lk}w}$ and
$\sigma_0^{uv^{i+1}w}=\sigma_0^{uv^{i+lk+1}w}$.
Hence for infinitely many natural numbers $n>m$, we have
$uv^nw \in L$ and $uv^{n+1}w \not\in L$.  
\end{proof}

The following lemma will be useful when
considering sets of geodesics in
subgroups.

\begin{lemma}\label{subalpha} 
Let $B$ be a subset of $A$.  Then $B^*$ is
a star-free language over $A$. 
\end{lemma}

\begin{proof}
The set $B^*$ can be written as
$$
B^* = ( \cup_{a \in A \setminus B} \emptyset^c a  \emptyset^c )^c,
$$
where $S^c$ denotes the complement $A^* \setminus S$ of a subset $S$ of
$A^*$.
This is a star-free expression for $B^*$.
\end{proof}

In Section~\ref{grprod} we will utilize another equivalent
characterization of star-free languages, shown in the
following proposition.

\begin{definition} We define a circuit in an automaton 
through a state $\sigma$ to
be {\em powered} if the circuit is labelled by a word $v^k$
for some $v \in A^+$
with $k>1$ 
such that $\sigma^v \neq \sigma^{v^2}$.
\end{definition}

\begin{proposition}\label{power_circuits}
A regular language $L$ is star-free if and only if the minimal automaton to
recognize it has no powered circuits.
\end{proposition}
\begin{proof}
Suppose that $v^k$ labels a powered
circuit in the minimal automaton $M_L$
for $L$ beginning at a state $\sigma$, with $k>1$ the least
natural number such that $\sigma^{v^k}=\sigma$.
Let $u$ label a route from the start state $\sigma_0$ to $\sigma$.
Since $M_L$ is minimal, the states $\sigma^u$ and $\sigma^{uv}$
must have different futures, so
there is a word $w$ such that one of $uw$ and $uvw$ is in $L$ and
the other is not.
Hence for all natural numbers $a$, we also have that
one of $uv^{ak}w$ and $uv^{ak+1}w$ is in $L$ and the other is not.
So, by Proposition~\ref{star-free1}, $L$ cannot be star-free.

Conversely, suppose that $L$ is not star-free. 
Let $m$ be the number of states of $M_L$.
Again using
Proposition~\ref{star-free1}, there exist words $u,v,w$ in $A^*$
such that $uv^nw \in L$ but $uv^{n+1}w \not \in L$ for some $n>m$.
Among the targets from the start state under
the transitions labelled $u, uv,
uv^2 ,\ldots, uv^{m}$ there must be at least one coincidence, and
hence there exist numbers $0 \le a<b \le m$ such that the targets
satisfy
$\sigma_0^{uv^a}=\sigma_0^{uv^b}$. Since $uv^nw$ and $uv^{n+1}w$
do not have the same target, neither do $uv^n$ and $uv^{n+1}$,
and hence we cannot have $b=a+1$.  Therefore
$v^{b-a}$ labels a powered circuit at $\sigma_0^{uv^a}$.
\end{proof}

\section{Small cancellation result}
\label{sc}

For a finite group presentation $\langle\, X \mid R\, \rangle$ in which
the elements of $R$ are cyclically reduced, we define
the {\em symmetrization} $R_*$ of $R$ to be the set of all cyclic
conjugates of all words in $R \cup R^{-1}$.

We define a {\em piece} to be a word over $X$ which is a prefix of two
(or more) distinct words in $R_*$.
Geometrically, a piece is a word labelling a common face
of two 2-cells, or regions, of a Dehn (or van Kampen) diagram
for this presentation.

Where $\lambda>0$,
we say that the presentation satisfies $C'(\lambda)$
if every piece has length less than $\lambda$ times the length of any
relator containing it; in other words, a word labelling a common
face between two diagram regions has length less than $\lambda$
times the length of either of the words labelling the boundaries
of the regions.

Where $q$ is a positive integer greater than 3,
we say that the presentation satisfies $T(q)$ if, for any
$h$ with $3 \le h < q$ and $r_1,\ldots ,r_h \in R_*$ with
$r_i \neq r_{i-1}^{-1}$ ($1<i\le h$) and $r_1 \ne r_h^{-1}$,
at least one of the words $r_1r_2$, $r_2r_3$, $\ldots$, $r_hr_1$
is freely reduced.  Geometrically, this implies that
each interior vertex of a Dehn diagram with degree greater than
2 must actually have degree at least $q$.
We refer the reader to~\cite{Strebel} for more details on small cancellation
conditions and diagrams.

In this section we prove the following theorem, utilizing results
from~\cite{Strebel} which show that sufficiently restrictive small 
cancellation conditions force Dehn diagrams to be very thin.

\begin{theorem} \label{thm:smallcanc}
Let $G = \langle\, X \mid R\, \rangle$ be a finite group
presentation such that the symmetrization $R_*$ of $R$ satisfies one of the
small cancellation conditions
$C'(1/6)$ (hypothesis A), or $C'(1/4)$ and $T(4)$ (hypothesis B).
Then the language of all words over $X$ that are geodesic in $G$
 is regular and star-free.
\end{theorem}
Note that it is already well known that these small cancellation 
conditions imply word hyperbolicity (see \cite{Strebel}), and hence
that the geodesic words form a regular set.

The remainder of this section is devoted to the proof of
Theorem~\ref{thm:smallcanc}.
Let $L$ be the language of all geodesic words for a presentation of a group $G$
satisfying the hypotheses of this theorem, and suppose to the
contrary that $L$ is not star-free.

For any word $y$ over $X \cup X^{-1}$, let $|y|$ denote
the length of the word $y$, and let $|y|_G$ denote
the length of the element of $G$ represented by $y$.

By Proposition~\ref{star-free1}
there exist words $u,v,w$ over $X$ such that there exist arbitrarily large $n$
with $uv^{n-1}w$ geodesic and $uv^nw$ not geodesic.
Choose $u$, $v$, and $w$ with this property such that $|u| + |w|$ is minimal.
Note that $u$, $v$, and $w$ are each nonempty.
Then there exist arbitrarily large $n$ with $uv^{n-1}w$ geodesic
and $uv^nw$ {\em minimally} non-geodesic (that is, all of its proper
subwords are geodesic), since otherwise $u',v,w$ or $u,v,w'$
would have the property in question, where $u'$ and $w'$ are respectively
the maximal proper suffix of $u$ and the maximal proper prefix of $w$.

Choose some $n$ for which $uv^nw$ is  minimally non-geodesic and
$uv^{n-1}w$ is geodesic.

Let $w = w'x$ with $x \in X \cup X^{-1}$, and let $t$ be the geodesic word
$uv^nw'$.
Then either
\begin{enumerate}
\item[(i)] $|uv^nw|_G = |t| - 1$, and there is a
geodesic word $t'$ ending in $x^{-1}$ with $|t'|= |t|$ and $t' =_G t$; or
\item[(ii)] $|uv^nw|_G = |t|$ and there is a geodesic word $t'$ with
$|t'|= |t|$ and $t' =_G tx$.
\end{enumerate}

In the first case, we let $T$ be the geodesic digon in the Cayley
graph $\Gamma$ of $G$ with edges labelled $t$ and $t'$,
and in the second case we let $T$ be the geodesic triangle
in $\Gamma$ with edges labelled $t$, $x$ and $t'$. In both cases, the
edges labelled $t$ and $t'$ start at the base point of $\Gamma$.

Notice that none of the internal vertices of the paths labelled
$t$ and $t'$ in $\Gamma$ can be equal
to each other, because such an equality would imply that $uv^nw$ 
has a proper suffix that is not geodesic, contrary to the choice of $n$.
Hence in any Dehn diagram with boundary
labels given by the words labelling the edges of $T$, the
vertices of the boundary paths labelled
$t$ and $t'$ also cannot be equal except at the endpoints.

It follows from the proof of Proposition 39 (ii) of~\cite{Strebel}
that there is a reduced Dehn diagram $\Delta$ with boundary
labels given by $T$, which has the form of one of
the two diagrams (corresponding to cases (i) and (ii) above) in Figure 10
of that proof, reproduced here in Figure 1.

\begin{center}
\begin{picture}(360,150)(0,-10)

\put(10,110){\line(2,1){10}}
\put(30,120){\line(2,1){10}}
\put(40,125){\line(1,0){10}}
\put(60,125){\line(1,0){10}}
\put(80,125){\line(1,0){10}}
\put(100,125){\line(1,0){10}}
\put(120,125){\line(1,0){10}}
\put(140,125){\vector(1,0){12}}
\put(150,128){$t$}
\put(152,125){\line(1,0){8}}
\put(170,125){\line(1,0){10}}
\put(190,125){\line(1,0){10}}
\put(210,125){\line(1,0){10}}
\put(230,125){\line(1,0){10}}
\put(250,125){\line(1,0){10}}
\put(260,125){\line(2,-1){10}}
\put(280,115){\line(2,-1){10}}

\put(10,100){\line(2,-1){10}}
\put(30,90){\line(2,-1){10}}
\put(40,85){\line(1,0){10}}
\put(60,85){\line(1,0){10}}
\put(80,85){\line(1,0){10}}
\put(100,85){\line(1,0){10}}
\put(120,85){\line(1,0){10}}
\put(140,85){\vector(1,0){12}}
\put(145,72){$t'$}
\put(152,85){\line(1,0){8}}
\put(170,85){\line(1,0){10}}
\put(190,85){\line(1,0){10}}
\put(210,85){\line(1,0){10}}
\put(230,85){\line(1,0){10}}
\put(250,85){\line(1,0){10}}
\put(260,85){\line(2,1){10}}
\put(280,95){\line(2,1){10}}

\put(10,105){\line(2,1){30}}
\put(10,105){\line(2,-1){30}}
\put(40,90){\line(0,1){30}}
\put(40,90){\line(1,0){25}}
\put(40,120){\line(1,0){25}}
\put(65,90){\line(0,1){30}}
\put(65,90){\line(1,0){25}}
\put(65,120){\line(1,0){25}}
\put(90,90){\line(0,1){30}}
\put(90,90){\line(1,0){25}}
\put(90,120){\line(1,0){25}}
\put(115,90){\line(0,1){30}}

\put(128,92){\circle*{2}}
\put(150,92){\circle*{2}}
\put(172,92){\circle*{2}}
\put(128,118){\circle*{2}}
\put(150,118){\circle*{2}}
\put(172,118){\circle*{2}}
\put(128,105){\circle*{2}}
\put(150,105){\circle*{2}}
\put(172,105){\circle*{2}}

\put(185,90){\line(0,1){30}}
\put(185,90){\line(1,0){25}}
\put(185,120){\line(1,0){25}}
\put(210,90){\line(0,1){30}}
\put(210,90){\line(1,0){25}}
\put(210,120){\line(1,0){25}}
\put(235,90){\line(0,1){30}}
\put(235,90){\line(1,0){25}}
\put(235,120){\line(1,0){25}}
\put(260,90){\line(0,1){30}}
\put(290,105){\vector(-2,-1){18}}
\put(272,96){\line(-2,-1){12}}
\put(260,120){\line(2,-1){30}}
\put(269,100){$x$}

\put(310,105){Case (i)}

\put(10,30){\line(2,1){10}}
\put(30,40){\line(2,1){10}}
\put(40,45){\line(1,0){10}}
\put(60,45){\line(1,0){10}}
\put(80,45){\line(1,0){10}}
\put(100,45){\line(1,0){10}}
\put(120,45){\line(1,0){10}}
\put(140,45){\vector(1,0){12}}
\put(150,48){$t$}
\put(152,45){\line(1,0){8}}
\put(170,45){\line(1,0){10}}
\put(190,45){\line(1,0){10}}
\put(210,45){\line(1,0){10}}
\put(230,45){\line(1,0){10}}
\put(250,45){\line(1,0){10}}
\put(260,45){\line(3,-1){10}}
\put(280,38.7){\line(3,-1){10}}

\put(10,20){\line(2,-1){10}}
\put(30,10){\line(2,-1){10}}
\put(40,5){\line(1,0){10}}
\put(60,5){\line(1,0){10}}
\put(80,5){\line(1,0){10}}
\put(100,5){\line(1,0){10}}
\put(120,5){\line(1,0){10}}
\put(140,5){\vector(1,0){12}}
\put(145,-8){$t'$}
\put(152,5){\line(1,0){8}}
\put(170,5){\line(1,0){10}}
\put(190,5){\line(1,0){10}}
\put(210,5){\line(1,0){10}}
\put(230,5){\line(1,0){10}}
\put(250,5){\line(1,0){10}}
\put(260,5){\line(3,1){10}}
\put(280,11.7){\line(3,1){10}}

\put(10,25){\line(2,1){30}}
\put(10,25){\line(2,-1){30}}
\put(40,10){\line(0,1){30}}
\put(40,10){\line(1,0){25}}
\put(40,40){\line(1,0){25}}
\put(65,10){\line(0,1){30}}
\put(65,10){\line(1,0){25}}
\put(65,40){\line(1,0){25}}
\put(90,10){\line(0,1){30}}
\put(90,10){\line(1,0){25}}
\put(90,40){\line(1,0){25}}
\put(115,10){\line(0,1){30}}

\put(128,12){\circle*{2}}
\put(150,12){\circle*{2}}
\put(172,12){\circle*{2}}
\put(128,38){\circle*{2}}
\put(150,38){\circle*{2}}
\put(172,38){\circle*{2}}
\put(128,25){\circle*{2}}
\put(150,25){\circle*{2}}
\put(172,25){\circle*{2}}

\put(185,10){\line(0,1){30}}
\put(185,10){\line(1,0){25}}
\put(185,40){\line(1,0){25}}
\put(210,10){\line(0,1){30}}
\put(210,10){\line(1,0){25}}
\put(210,40){\line(1,0){25}}
\put(235,10){\line(0,1){30}}
\put(235,10){\line(1,0){25}}
\put(235,40){\line(1,0){25}}
\put(260,10){\line(0,1){30}}
\put(260,10){\line(3,1){30}}
\put(260,40){\line(3,-1){30}}
\put(290,30){\vector(0,-1){7}}
\put(290,23){\line(0,-1){3}}
\put(294,22){$x$}

\put(310,25){Case (ii)}
\end{picture}
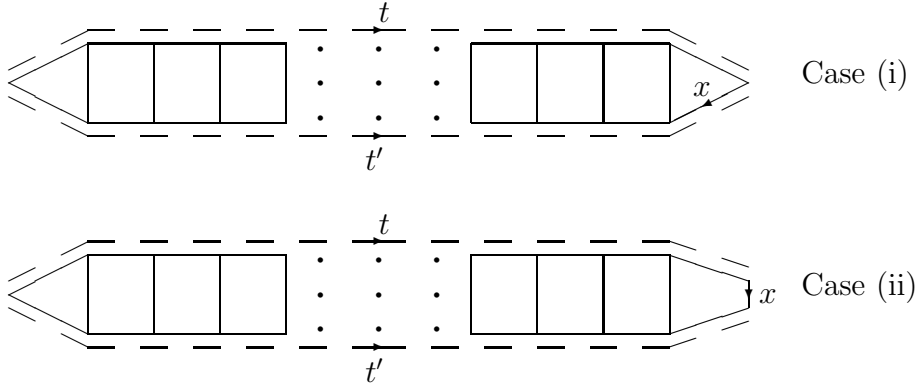
Figure 1.
\end{center}

A key feature is that all internal vertices in these Dehn diagrams 
have degree 2, and all external vertices have degree 2 or 3.
Notice that the short vertical paths joining vertices of $t$ to vertices
of  $t'$ in these diagrams are all pieces of the relators corresponding to the
two adjoining regions and hence have length less than $1/6$ or $1/4$ of those
relators under hypotheses A and B respectively.  

The boundary label of each region (except for
those regions containing the endpoints of $t$ and $t'$) has the form $r =
sps'p' \in R_*$, where $s$ and $s'$ are nonempty 
subwords of $t$ and $(t')^{-1}$,
respectively, and $p$ and $p'$ are nonempty pieces.
We consider $s$ to label the `top' of $r$, $s'^{-1}$ the bottom 
and $p$ and $p'^{-1}$
the right and left boundaries of $r$, respectively.
We shall use this top, bottom, right, and left convention throughout this
section when referring to any region or union
of consecutive regions of $\Delta$.

We shall establish the required contradiction to our assumption
that $L$ is not star-free by identifying a union  $\Psi_j$ of
consecutive regions of $\Delta$ such that the top of $\Psi_j$ is labelled by
a cyclic conjugate of $v$, the labels of the left and
right boundaries are identical, and the bottom label
is no shorter than $v$. 
In that case the diagram formed by deleting $\Psi_j$
from $\Delta$ and identifying its left and right boundaries
demonstrates that $uv^{n-1}w$ is not geodesic, contrary to our choice
of $u$, $v$, and $w$.

To construct the top boundary of our region $\Psi_j$, we start by finding a
subword $s$ of $v^n$ with particular properties.
We define a subword $s$ of $v^n$ to have the property ($\dagger$) if 
\begin{description}
\item[(i)] some occurrence of $s$ in $v^n$ labels 
the top of a region in $\Delta$ 
with boundary $r=sps'p'$ and
\item[(ii)] $|s| >|r|/3$ if hypothesis A holds, 
$|s|> |r|/4$ if hypothesis B holds.
\end{description}

\begin{lemma}
Provided that $n$ is large enough, $v^n$ has a 
subword $s$ satisfying ($\dagger$).
\end{lemma}
\begin {proof} Let $r=sps'p'$ label any region in $\Delta$. 
If $|s| \le |r|/3$ (resp. $|s| \le |r|/4$),
then since $|p|,|p'| < |r|/6$ (resp. $|p|,|p'| < |r|/4$), we have
$|s'| > |s|$. But, since $|t| = |t'|$, then provided that $n$ is large
enough, this cannot be true for all regions of $\Delta$
whose top label is a subword of $v^n$.
\end{proof}

From now on let $n$ be large enough for the lemma above to be applied.
In addition let $n$ be larger than $6\rho$, where 
$\rho$ is the length of the longest relator in the presentation.

For the remainder of this section, let $s$ be a specific choice
of a word which is as long as possible
subject to satisfying ($\dagger$).  
We shall call an occurrence of $s$ in $v^n$ {\em strictly internal} if it 
does not include any of the first $\rho$ or last $\rho$ letters of $v^n$.

\begin{lemma}\label{lem:internal}
Every strictly internal occurrence of $s$
in $v^n$ labels the top of a region in $\Delta$.
\end{lemma}

\begin{proof}
Suppose not, and consider a strictly internal occurrence of $s$ which
is not the top of a relator.

First note that $s$ cannot be a proper subword of the top $y$
of a region of $\Delta$. For if it were, 
since $s$ is strictly internal, the word $y$ would be a 
subword of $v^n$ satisfying ($\dagger$) and
have length longer than $s$.

Hence there are adjacent regions of $\Delta$ with labels
$r_1 = s_1p_1s'_1p'_1$ and $r_2 = s_2p_2s'_2p'_2$ where  $p'_2 = p_1^{-1}$
and the vertex at the end of the path labelled $s_1$ and at the beginning
of the path labelled $s_2$ is an internal vertex of the subpath of
$t$ labelled $s$.

Either (1) $s_1$ is a proper subword of $s$  
or (2) $s_2$ is a proper subword of $s$,
or (3) $s$ is a subword of $s_1s_2$ (see Figure 2).

\begin{center}
\begin{picture}(400,80)(-10,-12)
\put(0,10){\line(1,0){10}}
\put(0,65){\line(1,0){10}}
\put(10,10){\vector(0,1){30}}
\put(14,35){$p_1'$}
\put(10,40){\line(0,1){25}}
\put(10,65){\vector(1,0){38}}
\put(42,55){$s_1$}
\put(48,65){\line(1,0){32}}
\put(80,65){\vector(0,-1){28}}
\put(66,35){$p_1$}
\put(80,37){\line(0,-1){27}}
\put(80,10){\vector(-1,0){38}}
\put(42,-2){$s_1'$}
\put(42,10){\line(-1,0){32}}

\put(83,10){\vector(0,1){28}}
\put(87,35){$p_2'$}
\put(83,38){\line(0,1){27}}
\put(83,65){\vector(1,0){38}}
\put(115,55){$s_2$}
\put(121,65){\line(1,0){32}}
\put(153,65){\vector(0,-1){30}}
\put(139,35){$p_2$}
\put(153,35){\line(0,-1){25}}
\put(153,10){\vector(-1,0){38}}
\put(115,-2){$s_2'$}
\put(115,10){\line(-1,0){32}}
\put(153,10){\line(1,0){10}}
\put(153,65){\line(1,0){10}}

\put(-5,69){\vector(1,0){60}}
\put(51,72){$s$}
\put(55,69){\line(1,0){50}}
\put(83,-12){$(1)$}

\put(200,10){\line(1,0){10}}
\put(200,65){\line(1,0){10}}
\put(210,10){\vector(0,1){30}}
\put(214,35){$p_1'$}
\put(210,40){\line(0,1){25}}
\put(210,65){\vector(1,0){38}}
\put(242,55){$s_1$}
\put(248,65){\line(1,0){32}}
\put(280,65){\vector(0,-1){28}}
\put(266,35){$p_1$}
\put(280,37){\line(0,-1){27}}
\put(280,10){\vector(-1,0){38}}
\put(242,-2){$s_1'$}
\put(242,10){\line(-1,0){32}}

\put(283,10){\vector(0,1){28}}
\put(287,35){$p_2'$}
\put(283,38){\line(0,1){27}}
\put(283,65){\vector(1,0){38}}
\put(315,55){$s_2$}
\put(321,65){\line(1,0){32}}
\put(353,65){\vector(0,-1){30}}
\put(339,35){$p_2$}
\put(353,35){\line(0,-1){25}}
\put(353,10){\vector(-1,0){38}}
\put(315,-2){$s_2'$}
\put(315,10){\line(-1,0){32}}
\put(353,10){\line(1,0){10}}
\put(353,65){\line(1,0){10}}

\put(240,69){\vector(1,0){45}}
\put(281,72){$s$}
\put(285,69){\line(1,0){40}}
\put(283,-12){$(3)$}

\end{picture}
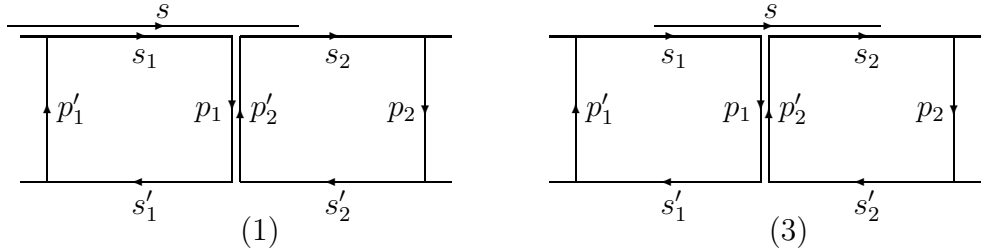
Figure 2.
\end{center}

Note that it is impossible for the $T(4)$ property 
to hold in any of these situations.
A violation is provided by the two relators shown, together with a relator
containing $s$. 
So we may assume that $C'(1/6)$ holds.

In case (1),  $s_1$ is a piece, and so are $p_1'$ and $p_1$
(since the $p_i$ and $p_i'$ are the labels of the vertical paths in $\Delta$),
and so each must have length less than $|r|/6$. But then the path
$p_1's_1p_1$ is shorter than $s_1'$, contradicting the fact that $t'$ is
a geodesic. Case (2) is dealt with similarly.

In case (3), $s \cap s_1$ and $s \cap s_2$ are pieces, and so must each have
length less than $|r|/6$, contradicting the condition that $|s|>|r|/3$.
\end{proof}

Now write $t = u v^\rho \tilde v v^\rho w'$, where
$\tilde v=v^{n-2\rho}$.  The first occurrence of $s$ in $\tilde v$
must start before the end of the first $v$; that is, $\tilde v$
has a prefix of the form $qs$ where $q$ is a proper prefix of $v$.
Since $n>6\rho$,
the suffix $v^{n-2\rho-3}$ of $\tilde v$ has length at least $4\rho-3 \ge \rho$,
and so is longer than $s$.
Hence this condition on $n$ ensures that
the words $vqs$ and $v^2qs$ are also prefixes of $\tilde v$.
Then all three of these occurrences
of $s$ are also strictly internal subwords in $v^n$.
By Lemma~\ref{lem:internal}, these three occurrences of $s$ are
each the top label of a region of $\Delta$, so these subwords must
be disjoint.  Hence $t$ must have a prefix of the form
$u v^\rho qsq'sq's$, where the three subwords labelled $s$ are
the three discussed above, and $sq'$ is a cyclic conjugate of $v$.

Let the three regions of $\Delta$
whose tops are labelled by these
$s$ subwords be called $\Phi_0$, $\Phi_2$ and $\Phi_4$.
If $q'$ is nonempty, let 
$\Phi_1$ and $\Phi_3$ be the regions 
attached immediately to the left of $\Phi_2$ and $\Phi_4$, respectively
(see Figure 3).
Since the $s$ subwords are strictly internal, these regions cannot
contain the endpoints of $t$ or $t'$.
Then for $i=0,\ldots, 4$ the region $\Phi_i$ has a boundary label 
of the form $r_i:= s_ip_is_i'p_i'$.  Note that $s_0=s_2=s_4=s$.

\begin{center}
\begin{picture}(400,90)(-25,-5)

\put(-25,69){\vector(1,0){60}}
\put(32,74){$q'$}
\put(35,69){\line(1,0){45}}

\put(-22,37){\circle*{2}}
\put(-6,37){\circle*{2}}
\put(-22,13){\circle*{2}}
\put(-6,13){\circle*{2}}
\put(-22,62){\circle*{2}}
\put(-6,62){\circle*{2}}

\put(0,10){\line(1,0){10}}
\put(0,65){\line(1,0){10}}
\put(10,10){\vector(0,1){30}}
\put(14,35){$p_1'$}
\put(10,40){\line(0,1){25}}
\put(10,65){\vector(1,0){38}}
\put(42,55){$s_1$}
\put(48,65){\line(1,0){32}}
\put(80,65){\vector(0,-1){28}}
\put(66,35){$p_1$}
\put(80,37){\line(0,-1){27}}
\put(80,10){\vector(-1,0){38}}
\put(42,-2){$s_1'$}
\put(42,10){\line(-1,0){32}}
\put(42,34){$\Phi_1$}

\put(83,10){\vector(0,1){28}}
\put(87,35){$p_2'$}
\put(83,38){\line(0,1){27}}
\put(83,65){\vector(1,0){38}}
\put(100,70){$s_2=s$}
\put(121,65){\line(1,0){32}}
\put(153,65){\vector(0,-1){30}}
\put(139,35){$p_2$}
\put(153,35){\line(0,-1){25}}
\put(153,10){\vector(-1,0){38}}
\put(115,-2){$s_2'$}
\put(115,10){\line(-1,0){32}}
\put(115,34){$\Phi_2$}
\put(153,10){\line(1,0){10}}
\put(153,65){\line(1,0){10}}

\put(153,69){\vector(1,0){66}}
\put(212,74){$q'$}
\put(219,69){\line(1,0){61}}

\put(168,37){\circle*{2}}
\put(182,37){\circle*{2}}
\put(196,37){\circle*{2}}
\put(168,13){\circle*{2}}
\put(182,13){\circle*{2}}
\put(196,13){\circle*{2}}
\put(168,62){\circle*{2}}
\put(182,62){\circle*{2}}
\put(196,62){\circle*{2}}

\put(200,10){\line(1,0){10}}
\put(200,65){\line(1,0){10}}
\put(210,10){\vector(0,1){30}}
\put(214,35){$p_3'$}
\put(210,40){\line(0,1){25}}
\put(210,65){\vector(1,0){38}}
\put(242,55){$s_3$}
\put(248,65){\line(1,0){32}}
\put(280,65){\vector(0,-1){28}}
\put(266,35){$p_3$}
\put(280,37){\line(0,-1){27}}
\put(280,10){\vector(-1,0){38}}
\put(242,-2){$s_3'$}
\put(242,10){\line(-1,0){32}}
\put(242,34){$\Phi_3$}

\put(283,10){\vector(0,1){28}}
\put(287,35){$p_4'$}
\put(283,38){\line(0,1){27}}
\put(283,65){\vector(1,0){38}}
\put(300,70){$s_4=s$}
\put(321,65){\line(1,0){32}}
\put(353,65){\vector(0,-1){30}}
\put(339,35){$p_4$}
\put(353,35){\line(0,-1){25}}
\put(353,10){\vector(-1,0){38}}
\put(315,-2){$s_4'$}
\put(315,10){\line(-1,0){32}}
\put(315,34){$\Phi_4$}
\put(353,10){\line(1,0){10}}
\put(353,65){\line(1,0){10}}

\end{picture}
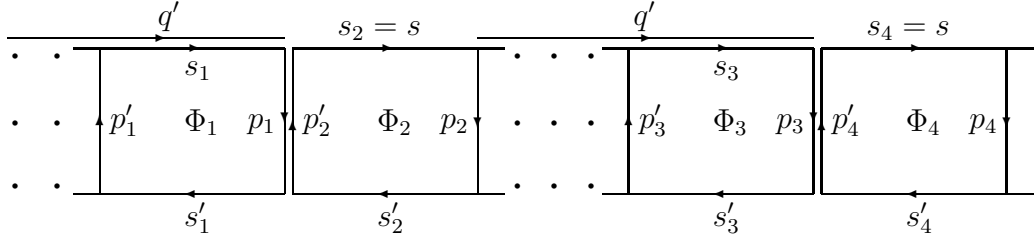
Figure 3.
\end{center}

\begin{lemma} We have $r_2=r_4$ and $p_2' = p_4'$.
\end{lemma}
\begin{proof}
Since $s$  satisfies ($\dagger$) and so is too long to be a piece, it follows
immediately that $r_2=r_4$.

We shall show that $p_2' = p_4'$ in the case when $q'$ is not the empty word;
the proof in the case when $q'$ is empty is similar.

First we shall prove that $s_1=s_3$.

Under hypothesis A, since $p_1$ and $p_1'$ are pieces and
$t'$ is geodesic, we have $|p_1|,|p_1'| < |r_1|/6$ and
$|s_1'| < |r_1|/2$, so therefore $|s_1| > |r_1|/6$.
Thus $s_1$ is too long to be a piece, so the relator $r_1$ is
is the only element of $R_*$ with prefix $s_1$.
Let $y$ be the longest suffix of $q'$ which is
a subword of an element of $R_*$.  Then 
$s_1$ is a suffix of $y$, i.e.~$y=zs_1$, and 
again $y$ is too long to be a piece, so the only relators
in $R_*$ containing $y$ are cyclic conjugates of $r_1$.
Hence the (possibly empty) word $z$ 
is a suffix of $r_1$.
If $z$ were nonempty, then 
the last letter $a$ of $z$ would
be a suffix of $p_1'$, and
the relator 
$\tilde r=\tilde s \tilde p \tilde s' \tilde p'$ labelling
the boundary of the region $\tilde \Phi$ immediately to the
left of $\Phi_1$ would have top label $\tilde s$ 
ending with $a$.  But the first letter
of $\tilde p=p_1'^{-1}$ would then be $a^{-1}$, contradicting
the fact that the relator $\tilde r \in R_*$ must
be reduced.  Thus $s_1 $ must be the longest suffix
of $q'$ which is
a subword of an element of $R_*$ in this case.

Under hypothesis B, 
again let $y$ be the longest suffix of $q'$ which is
a subword of an element of $R_*$, and write $y=zs_1$.
If $z$ were nonempty, then the regions $\tilde \Phi$,
$\Phi_1$, and a third region with boundary label
containing $y$ glued along this word to the other two
regions, would violate the $T(4)$ condition.  
Thus again $s_1 $ must be the longest suffix
of $q'$ which is
a subword of an element of $R_*$.

The same argument under both hypotheses
shows that $s_3$ is also the longest suffix
of $q'$ which is
a subword of an element of $R_*$, so $s_1=s_3$.

Now suppose that $p_2' \neq p_4'$. Then one is a subword of the other, since
both are suffixes of $r_2=r_4$. Suppose (without loss of generality)
that $p_4'=\tau^{-1}p_2'$; then
$p_3=p_1\tau$. If $r_1=r_3$ then $\tau$ is a prefix of 
$s_1'p_1's_1p_1$ and
$\tau^{-1}$ is a suffix of $p_2's_2p_2s_2'$, and so
$t'$ is not freely reduced. If $r_1\neq r_3$,
then $s_1p_1=s_3p_1$ is a piece, and under 
either $C'(1/6)$ or $C'(1/4)$ must
have length less than $|r_1|/4$. 
Then $p_1's_1p_1$ has length less than $|r_1|/2$,
and $s_1'$ cannot be geodesic.
\end{proof}

Now let $\Psi_1$ be the union of the regions $\Phi_2$ and all regions
of $\Delta$ to the right of $\Phi_2$
up to but not including $\Phi_4$. Then the boundary
label of $\Psi_1$ is $sq' p_4'^{-1} \sigma'_1 p_2'$ where
$sq'$ is a cyclic conjugate of $v$, 
$p_4'  = p_2'$, and $\sigma_1'$ is a subword of $(t')^{-1}$
(see Figure 4 with $q' \neq 1$).

\begin{center}
\begin{picture}(370,110)(15,-10)

\put(23,10){\line(1,0){10}}
\put(23,65){\line(1,0){10}}
\put(33,10){\vector(0,1){28}}
\put(37,35){$p_2'$}
\put(33,38){\line(0,1){27}}
\put(33,65){\vector(1,0){38}}
\put(50,70){$s_2=s$}
\put(71,65){\line(1,0){32}}
\put(103,65){\vector(0,-1){30}}
\put(103,35){\line(0,-1){25}}
\put(103,10){\vector(-1,0){38}}
\put(65,10){\line(-1,0){32}}
\put(65,34){$\Phi_2$}
\put(103,10){\line(1,0){10}}
\put(103,65){\line(1,0){10}}

\put(103,69){\vector(1,0){90}}
\put(180,74){$q'$}
\put(193,69){\line(1,0){87}}

\put(132,37){\circle*{2}}
\put(162,37){\circle*{2}}
\put(192,37){\circle*{2}}
\put(132,13){\circle*{2}}
\put(162,13){\circle*{2}}
\put(192,13){\circle*{2}}
\put(132,62){\circle*{2}}
\put(162,62){\circle*{2}}
\put(192,62){\circle*{2}}

\put(200,10){\line(1,0){10}}
\put(200,65){\line(1,0){10}}
\put(210,10){\vector(0,1){30}}
\put(210,40){\line(0,1){25}}
\put(210,65){\vector(1,0){38}}
\put(248,65){\line(1,0){32}}
\put(280,65){\vector(0,-1){28}}
\put(266,35){$p_3$}
\put(280,37){\line(0,-1){27}}
\put(280,10){\vector(-1,0){38}}
\put(242,10){\line(-1,0){32}}
\put(242,34){$\Phi_3$}

\put(283,10){\vector(0,1){28}}
\put(287,35){$p_4'$}
\put(283,38){\line(0,1){27}}
\put(283,65){\vector(1,0){38}}
\put(302,70){$s_4=s$}
\put(321,65){\line(1,0){32}}
\put(353,65){\vector(0,-1){30}}
\put(353,35){\line(0,-1){25}}
\put(353,10){\vector(-1,0){38}}
\put(315,10){\line(-1,0){32}}
\put(315,34){$\Phi_4$}
\put(353,10){\line(1,0){10}}
\put(353,65){\line(1,0){10}}

\put(280,6){\vector(-1,0){130}}
\put(140,-5){$\sigma_1'$}
\put(150,6){\line(-1,0){117}}

\put(159,37){\oval(280,105)}
\put(165,-10){$\Psi_1$}

\end{picture}
\end{center}
\begin{center}
Figure 4.
\end{center}

Provided that $n$ is large enough, we can define similar unions
of regions $\Psi_2$, $\Psi_3$, $\ldots$, $\Psi_i$, $\ldots$ of $\Delta$,
where $\Psi_i$ is immediately to
the left of $\Psi_{i+1}$, and $\Psi_i$ has boundary label
$sq'p_4'^{-1} \sigma'_i p'_2$, where each $\sigma'_i$ is a subword of
$(t')^{-1}$.

We have not attempted to show that the words $\sigma'_i$ are equal for all
$i$, but since $|t| = |t'|$ we cannot have $|\sigma'_i| < |v|$ for all $i$,
and so there exists a $j$ with  $|\sigma'_j| \ge |v|$. Now, if we
remove $\Psi_j$ from $\Delta$, the effect is to replace $t$ by $uv^{n-1}w'$
and $t'$ by a word of length at most $|t| - |v|$, so we obtain a
diagram that shows that $uv^{n-1}w$ is not geodesic, contrary
to assumption.

\section{Dependence on generating set}
\label{dependence}

In this section, we give two examples of groups for which the set of 
geodesics is star-free with respect to one generating set, but not 
with respect to another. 

Free groups have star-free geodesics with respect to
their free generators, but the following example shows that this is not
necessarily the case with an arbitrary generating set.  Let

\label{example}
 \[ G := \langle\, a,b,c,d,r,s \mid ba^2d = rcs,\, bd=s\, \rangle. \]

The two relations can be written as $r = b a^2 b^{-1} c^{-1}$, $s = bd$, and
so they can be used to eliminate $r$ and $s$ from any word representing an
element of $G$.  Hence $G$ is free on $a,b,c$ and $d$.

Since $ba^{2k}d =_G (ba^2b^{-1})^kbd =_G (rc)^ks$ for all $k \ge 0$,
we have that $ba^{2k}d$ is not a geodesic word.
We shall now show that $ba^{2k+1}d$ is a geodesic word for all $k \ge 0$,
which implies, by Proposition~\ref{star-free1},
that the set of geodesics for the group defined by this
presentation is not star-free.

Let $w = ba^{2k+1}d$ for some $k \ge 0$.  Suppose to the
contrary that $w$ is not geodesic.  Then there is a word  $x$  in
$a,b,c,d,r,s$ and their inverses satisfying $l(x) < l(w)$ which 
freely reduces to
$w$ after we make the above substitutions to eliminate $r$ and $s$.
Let $p_a$ denote the number of occurrences of $a$
in $x$, let $n_a$ denote the number of occurrences of $a^{-1}$
in $x$, and similarly for the other five generators of $G$.
Since the exponent sum of the $a$'s in $w$ is $2k+1$, 
and the only letters of $x$ that contribute powers of 
$a$ after the substitution are powers of $a$ and $r$, we have
$p_a-n_a+2p_r-2n_r=2k+1$.  Computing the
exponent sums of the $b$'s, $c$'s, and $d$'s in $w$ in the same way,
we obtain $p_b-n_b+p_s-n_s=1$, $p_c-n_c-p_r+n_r=0$,
and $p_d-n_d+p_s-n_s=1$, respectively.
Then 
\begin{eqnarray*}
2k+3 = l(w) 
   &  > & l(x) \\
    &=&   p_a+n_a+p_b+n_b+p_c+n_c+p_d+n_d+p_r+n_r+p_s+n_s\\
     &=& 2n_a+2n_b+2n_c+p_d+n_d+2n_r+2n_s+2k+2.\\
\end{eqnarray*}
Therefore $n_a=n_b=n_c=n_d=n_r=n_s=p_d=0$, so
the word $x$ contains no occurrences of
inverses of the generators of $G$, and also no occurrences of the
generator $d$.
As a consequence, then $p_a+2p_r=2k+1$, $p_s=1$, $p_b=0$, 
and $p_r=p_c$.  Hence the letter $b$ also
does not occur in $x$, the letter $s$ occurs once,
and since $2k+1$ is odd, we have $p_a>0$ so $x$ contains
at least one $a$.
Since the  $d$ is at the right hand end of $w$, the  $s$  must be at the
right hand end of  $x$.  Hence $x$  has the form  $y a z s$,
where $y$ and $z$ are words over $a,c,r$.
After making the substitutions for $r$ and $s$ in $x$,
all of the powers of $a$ in the word are positive, so the $a$ in the  
expression $y a z s$ 
for $x$ will not be cancelled after further free reduction,
but the exponent sum of the $b$'s to the left of this $a$ is zero. 
Thus the resulting word cannot freely reduce to $w$, giving the
required contradiction.

Our second example is the three-strand braid group $B_3$, which has
a presentation
$\langle\,a,b\mid bab=aba\,\rangle$. The geodesics for this group
on generators $\{a^{\pm 1},b^{\pm 1}\}$ are described in~\cite{Sabalka}. 
A reduced
word is geodesic if it does not contain both one of $\{ab,ba\}$ and also one
of $\{a^{-1}b^{-1},b^{-1}a^{-1}\}$ as subwords, and it does not
contain both $aba$ and also one of $\{a^{-1}, b^{-1}\}$ as subwords, and it
does not contain both $a^{-1}b^{-1}a^{-1}$ and also one of $\{a,b\}$ 
as subwords.  Hence the language of geodesics can be expressed as
\begin{eqnarray*}
&&[(\emptyset^c ab \emptyset^c \cup \emptyset^c ba \emptyset^c)
 \cap 
 (\emptyset^c a^{-1}b^{-1} \emptyset^c \cup 
   \emptyset^c b^{-1}a^{-1} \emptyset^c)]^c\\
&&\quad \cap \quad 
[(\emptyset^c aba \emptyset^c)
 \cap 
 (\emptyset^c a^{-1} \emptyset^c \cup 
   \emptyset^c b^{-1} \emptyset^c)]^c\\
&&\quad \cap\quad  
[(\emptyset^c a^{-1}b^{-1}a^{-1} \emptyset^c)
 \cap 
 (\emptyset^c a \emptyset^c \cup 
   \emptyset^c b \emptyset^c)]^c,\\
\end{eqnarray*}
which is a star-free
regular language.

In~\cite{CharneyMeier}, Charney and Meier prove that Garside groups have
regular geodesics with respect to the generating set consisting of the
divisors of the Garside element and their inverses. 
The class of Garside groups includes
and generalizes the class of Artin groups of finite type, which itself includes
the braid groups.
The set of divisors of the Garside element in the three-strand
braid group $B_3$ is
$\{a,b,ab,ba,aba\}$, and an automaton accepting the geodesics in the
positive monoid of this example is calculated explicitly in
Example 3.5 of~\cite{CharneyMeier}. We see from this that
$(ba)(aba)^n(a)$ is a geodesic for $n$ even, but not for $n$ odd.
So, by Proposition~\ref{star-free1},
the language of geodesics for $B_3$ with this second generating
set is not star-free.

\section{Direct products, free products and graph products}\label{grprod}

We start by proving the straightforward result that the class of groups with
star-free sets of geodesics is closed under taking direct products.

\begin{lemma}\label{dp}
If the languages $L_1,L_2$ of all geodesics of
$(G_1,X_1)$ and $(G_2,X_2)$ are star-free then so is the language $L$
of all geodesics of
$(G_1 \times G_2, X_1 \cup X_2)$.
\end{lemma}

\begin{proof}
The language $L$ is the set of words over $X_1 \cup X_2$
which project onto words in each
of $L_1$ and $L_2$ if we map in turn the elements 
of $X_2,X_1$ to the empty strings.
We show that $L$ can be described as the intersection of two star-free
languages, $L_1'$ and $L_2'$.

The language $L_1'$ is defined by wrapping arbitrary strings in $X_2$ around
the elements of $X_1$ for each string in $L_1$; that is,
$$L_1' := \{w_0x_1w_1 \cdots x_nw_n ~|~ w_i \in X_2^* {\rm \ and\ } 
 x_1 \cdots x_n \in L_1\}.
$$
Thus a regular
expression for $L_1'$ is found 
by replacing each element $x$ of $X_1$ in a
star-free (regular) expression for $L_1$ by $X_2^*xX_2^*$. 
Note that Lemma~\ref{subalpha} shows that $X_2^*$ is a star-free
language over $X_1 \cup X_2$.
Since star-free languages
are closed under concatenation, this regular expression
for $L_1'$ shows that $L_1'$ is star-free.
The language $L_2'$ is defined similarly.
\end{proof}

It is straightforward to prove that an analogue of the above result also
holds for free products. In fact in the next theorem
we prove a more general result, which
includes both direct and free products as special cases.

\begin{definition}
Let $\Gamma$ be a finite undirected graph with $n$ vertices labelled by
finitely generated groups $G_i$. Then the {\em graph product} $\Pi_\Gamma G$
of the groups $G_i$ with respect to $\Gamma$ is defined to be the group
generated by $G_1,\ldots,G_n$ modulo relations implying that elements of
$G_i$ and $G_j$ commute if there is an edge in $\Gamma$ connecting the
vertices labelled by $G_i$ and $G_j$.
\end{definition}

Thus if $\Gamma$ is either 
a graph with no edges or a complete graph, then the
graph product is the free or the direct product of the 
groups $G_i$, respectively.

The word problem for graph products is studied in detail
in~\cite{Green} and~\cite{HermillerMeier}. 
If we use a generating set for $\Pi_\Gamma G$ that
consists of the union of generating sets of the vertex groups $G_i$, then it
turns out that a word $w$ in the generators is non-geodesic if and only if pairs
of adjacent generators in $w$ that lie in commuting pairs of vertex
groups can be swapped around so as to produce a non-geodesic
subword lying in one of the vertex groups. This property is used in~\cite{LMW}
to prove that the geodesics form a regular set if and only if the geodesics of
the vertex groups all form regular sets. We adapt this proof to show that the
same is true with `regular' replaced by `star-free'. 

\begin{theorem}
Let $G_1,\ldots G_n$ be the vertex groups of a graph product
$G := \Pi_\Gamma G$,
let $A_1,\ldots A_n$ be finite inverse-closed sets of generators for $G_i$,
and let $L_1,\ldots L_n$ be the languages of all geodesics
in $G_1,\ldots G_n$ over $A_1,\ldots A_n$, respectively.
Let $A := \cup _{i=1}^n A_i$ and let $L$ be the language
of all geodesics in $G$ over $A$.
The languages $L_1,\ldots L_n$ are all star-free if and only if 
the language $L$ is star-free.
\end{theorem}

\begin{proof}
First, suppose that $L$ is star-free.  For each vertex index $i$,
a word over $A_i^*$ which is
geodesic as an element of $G_i$ is also geodesic as an element of $G$,
and the language $L_i$ is the intersection of the 
star-free language $L$ with $A_i^*$.
Lemma~\ref{subalpha} shows that the set $A_i^*$ is star-free,
so $L_i$ is also star-free.

Conversely, suppose that each language $L_i$ is star-free.
Let $F_i$ denote the minimal finite state automaton over $A_i$ that
accepts $L_i$.  Since any prefix of a geodesic is also
geodesic, the language $L_i$ is prefix-closed, and
therefore the automaton $F_i$ has a single fail 
state, and all other states are accept states.

Following the proof in~\cite{LMW}, for each $i$ define a finite state
automaton $\hat{F_i}$ over $A$ by adding arrows for the
generators in $A \setminus A_i$ to the automaton $F_i$ as follows.
For each $a \in A \setminus A_i$ which commutes with $G_i$, a loop labelled
$a$ is added at every state of $F_i$ (including the fail state). For each
$b \in A \setminus A_i$ which does not commute with $G_i$, an arrow labelled
$b$ is added to join each accept state of $F_i$ to the  start state, 
and a loop labelled $b$ is added at the fail state.
Completing the construction in~\cite{LMW}, 
an automaton $F$ is built to accept the intersection of
the languages of the automata $\hat{F_i}$, and the authors 
show that the language accepted by $F$
is exactly the language $L$ of geodesics of $\Pi_\Gamma G$ over $A$.

Note that if two states of $\hat{F_i}$ have the same future, 
then these states have the
same future under the restricted alphabet $A_i$.
Thus minimality of the automaton
$F_i$ implies that $\hat{F_i}$ is also minimal.

Since $L_i$ is star-free,
Proposition~\ref{power_circuits} says that $F_i$
has no powered circuits.  Then the finite state automaton
$\hat{F_i}$ has no powered circuit labelled by a word in $A_i^+$.
Suppose that $\hat{F_i}$ has a powered circuit over $A^+$ and
let $v$ be a least length word such that a power of $v$
labels a powered circuit in $\hat{F_i}$.  
Let $\sigma$ be the beginning state of this circuit,
and let $k>1$ be the least natural number
such that $\sigma^{v^k}=\sigma$.  The states in this circuit
must be accept states.  
If $v$ contains a letter 
$a \in A \setminus A_i$ which commutes with $G_i$, then
the word $v$ with the letter $a$ removed also labels a
powered circuit at $\sigma$, contradicting the choice of
$v$ with least length.  Then we can write $v=v_1bv_2$ such that
$v_1,v_2 \in A^*$ and the letter
$b \in A \setminus A_i$ does not commute with $G_i$.
Hence $\sigma^{v_1b}$ is the start state $\sigma_{i,0}$
of $\hat{F_i}$ and the word $v_2v_1b$
labels a circuit at $\sigma_{i,0}$, and so $v$ labels a
circuit at $\sigma$, contradicting the condition $k>1$.
Therefore $\hat{F_i}$ also has no powered circuits.
Applying Proposition~\ref{power_circuits} again,
then the language accepted by $\hat{F_i}$ is star-free.

Finally, the language $L$ is the intersection of the star-free languages
accepted by the $\hat{F_i}$, so
therefore $L$ is also star-free.
\end{proof}

\section{Virtually abelian groups}
\label{va}

In this section, we prove that every finitely generated abelian group has a
star-free set of geodesics with respect to any finite generating set, 
whereas every
finitely generated virtually abelian group $G$ has some finite generating set
with respect to which $G$ has star-free geodesics.

As a special case of these results, note that we can see quickly
that the geodesic language of $\Z^n$ for the standard (inverse
closed) generating set is
star-free, either via Lemma~\ref{dp} or as follows.  
The minimal automaton for this language has
states corresponding to subsets of the generators that do not
contain inverse pairs, together with a fail state.  At 
a state given by a subset $S$, the transition corresponding
to a generator $a$ will go either to $S$ itself if $S$ contains
$a$, to $S \cup \{a\}$ if $S$ does not contain $a$ or $a^{-1}$,
and to the fail state if $S$ contains $a^{-1}$.  
Then the transition monoid of the minimal automaton
for this language of geodesics, i.e.~the syntactic
monoid, is both abelian and generated by
idempotents, and hence every element is an idempotent and the monoid is
aperiodic.  Sch\"utzenberger's Theorem then says that
this language is star-free.

Our arguments in this section make use of a condition
that is more restrictive than the star-free property, and
indeed more restrictive than the piecewise testable property, 
which we shall call {\em piecewise excluding}.
A language $L$ over $A$ is said to be piecewise excluding if there
is a finite set of strings $W \subset A^*$ with the property that a word
$w \in A^*$ lies in $L$ if and only if $w$ does not contain any of the
strings in $W$ as a not necessarily consecutive substring.
In other words,
$$L = ( \cup_{i=1}^n \{A^* a_{i1} A^* a_{i_2} A^* \cdots
  A^* a_{il_i} A^* \})^c,$$
where $W = \{a_1,\ldots,a_n\}$ and $a_i = a_{i1}a_{i_2}\cdots a_{il_i}$
for $1 \le i \le n$. It follows directly from the above expression and the
definition of piecewise testable languages given in Section~\ref{intro}
that piecewise excluding languages are piecewise testable, and hence star-free.

We also need the following technical result, in which the set $\N$ of
natural numbers includes 0.

\begin{lemma}\label{finmin}
Define the ordering $\preceq$ on $\N^r$
by $(m_1,\ldots,m_r) \preceq (n_1,\ldots,n_r)$ if and only if $m_i \le n_i$ for
$1 \le i \le r$. Then any subset of $\N^r$ has only finitely many elements
that are minimal under $\preceq$.
\end{lemma}
\begin{proof} See~\cite[Lemma 4.3.2]{ECHLPT} or the second paragraph
of the proof of~\cite[Proposition 4.4]{NeumannShapiro}.
\end{proof}

\begin{proposition}
\label{abelian}
If $G$ is a finitely generated abelian group, then the set of all geodesic
words for any finite monoid generating set of $G$ is a piecewise
excluding (and hence a piecewise testable and a star-free) language.
\end{proposition}
\begin{proof}
Let $A = \{x_1,\ldots,x_r\}$ be a finite monoid generating set for $G$,
and let $$U := \{\, (n_1,\ldots,n_r) \in \N^r \mid
x_1^{n_1}x_2^{n_2}\cdots x_r^{n_r}\  \hbox{is non-geodesic}\,\}.$$
By Lemma~\ref{finmin}, the subset $V \subseteq U$, consisting of those elements
of $U$ that are minimal under $\preceq$, is finite.
Let $W$ be the set of all permutations of all words
$x_1^{n_1}x_2^{n_2}\cdots x_r^{n_r}$ with $(n_1,\ldots,n_r) \in V$.
Since $G$ is abelian, whether or not a word over $A$ is geodesic is
not changed by permuting the generators in the word.
Hence, by definition of minimality under $\preceq$, any non-geodesic word
over $A$ contains a word in $W$ as a not necessarily consecutive substring.
Conversely, since the words in $W$ are themselves non-geodesic, any word over
$A$ containing one of them as a not necessarily consecutive substring is also
non-geodesic.  So the set of geodesic words over $A^*$ is piecewise
excluding, as claimed.
\end{proof}

\begin{proposition}
Any finitely generated virtually abelian group
has an inverse-closed generating set
with respect to which the set of all geodesics is a piecewise testable
(and hence a star-free) language.
\end{proposition}
\begin{proof}
The analogous result for `regular' rather than `piecewise testable' is proved
in Propositions 4.1 and 4.4 of~\cite{NeumannShapiro}.
We extend that proof.

Let $N$ be a finite index normal abelian subgroup in $G$. We choose a
finite generating set of $G$ of the form $Z = X \cup Y$ with the following
properties.
\begin{enumerate}
\item\label{p1}
$X \subset N$ and $Y \subset G \setminus N$.
\item\label{p2}
Both $X$ and $Y$ are closed under the taking of inverses.
\item\label{p3}
$X$ is closed under conjugation by the elements of $Y$.
\item\label{p4}
$Y$ contains at least one representative of each nontrivial
coset of $N$ in $G$.
\item\label{p5}
For any equation $w =_G xy$ with $w$ a word of length at most 3 over $Y$,
$y \in Y \cup\{1\}$ and $x \in N$, we have $x \in X$.
\end{enumerate}

We must first show that such generating sets exist. To see this, start with
any finite generating set $Z$ of $G$ and let $X := Z \cap N$,
$Y := Z \setminus X$. Adjoin finitely many new generators to ensure that
Property~\ref{p4} holds and that $Y$ is closed under taking inverses. 
Since $Y$ is finite, there are only finitely many possible words $w$ in
Property~\ref{p5}, and we can adjoin finitely many new generators in $N$ to
ensure that Property~\ref{p5} holds. Now adjoin inverses of elements of $X$
to get Property~\ref{p2}. Since $N$ is abelian and $|G:N|$ is finite,
elements of $N$ have only finitely many conjugates in $G$, and so we can
adjoin finitely many conjugates of elements of $X$ to $Z$ to get
Property~\ref{p3}, after which $X$ will still be closed under inversion.
The five properties will then all hold.

Now let $L$ be the set of all geodesic words over $Z$. For $i \ge 0$, let
$Z_i$ be the set of all words $z_1\cdots z_m \in Z^*$ for which precisely $i$
of the symbols $z_j$ lie in $Y$, and let $L_i := L \cap Z_i$.
Let $\widetilde Z_i$ be the set of words in $Z^*$ containing at
least $i$ letters of $Y$.  Then 
$\widetilde Z_i :=\cup_{y_1,...,y_i \in Y} \{Z^*y_1Z^* \cdots Z^*y_iZ^*\}$
is a piecewise testable language.  The set $Z_i$ equals the intersection
$\widetilde Z_i \cap (\widetilde Z_{i+1})^c$, 
and so is also piecewise testable.

First we shall show that $L= L_0 \cup L_1 \cup L_2$.
Property~\ref{p3} implies that any word
in $L$ is equal in $G$ to a word of the same
length involving the same elements of $Y$, but with all of those elements
appearing at the right hand end of the word. By Property~\ref{p5},
any word over $Y$ of length three or more is non-geodesic.
So $L$ is contained in (and hence equal to) $L_0 \cup L_1 \cup L_2$.

Now Property~\ref{p1} ensures that an element of $L_1$ cannot
represent an element of $N$, and the same is true for $L_2$, because
Property~\ref{p5} implies that a word in $Z_2$ that represents an
element of $N$ cannot be geodesic. So $L_0$ is equal to the set of all
geodesic words in $N$ over $X$, and the set $X$ generates
the subgroup $N$.  Then $L_0$ is piecewise testable by
Proposition~\ref{abelian}.

Next we show that $L_1$ is piecewise testable.
For a fixed $y \in Y$,
applying the same argument as in the proof of
Proposition~\ref{abelian} to $X^*y$, we can show
that there is a finite set $W_y$ of words over $X$ with the property that
a word in $X^*y$ is non-geodesic if and only if it contains one of the
words in $W_y$ as a not necessarily consecutive substring.

For each $y \in Y$ and $x \in X$, denote the generator in $X$ equal in $G$ to
$y^{-1}xy$ (which exists by Property~\ref{p3}) by $x^y$.
Define $W$ to be the (finite) set of words over $Z$ of the form
$x_1\cdots  x_t y x_{t+1}^y \cdots x_s^y$, where $y \in Y$ and
$x_1\cdots  x_t x_{t+1}\cdots x_s \in W_y$.
Then a word in $Z_1$ is non-geodesic if and only if it contains a word in $W$
as a not necessarily consecutive substring.  So $L_1 = P \cap Z_1$, where $P$
is a piecewise excluding language. Since $Z_1$ is piecewise testable,
this shows that $L_1$ is piecewise testable.

The proof that $L_2$ is piecewise testable is similar and is left to the reader.
So $L = L_0 \cup L_1 \cup L_2$ is piecewise testable.
\end{proof}

\end{document}